\newif\ifarxiv
\arxivtrue  
\ifarxiv
\documentclass{article}
\usepackage{amsthm}
\newtheorem{theorem}{Theorem}
\newtheorem{remark}{Remark}
\usepackage{authblk}
\usepackage[margin=1in]{geometry}
\else
\documentclass[runningheads]{LNCS/llncs}
\fi
\usepackage[T1]{fontenc}
\usepackage{graphicx}
\usepackage{xcolor}
\usepackage{graphicx} 
\usepackage{subfig}
\usepackage{subcaption}
\usepackage{amsmath}
\usepackage{amsfonts}
\usepackage{float}
\usepackage{algorithm}
\usepackage{algorithmicx}
\usepackage{algpseudocode}
\usepackage{mathrsfs}
\usepackage{comment}
\usepackage{url}
\usepackage{cleveref}
\usepackage{tikz,tikz-cd,tikz-3dplot}
\usetikzlibrary{arrows.meta, positioning, calc, positioning, backgrounds,intersections, patterns}
\usepackage{verbatim}
\usepackage{cite}
\usepackage{wrapfig}
\usepackage{pifont}
\usepackage{enumitem}

\newcommand{\drawTriangle}[4][white]{
    \fill[#1] (#2) -- (#3) -- (#4) -- cycle;  
    \draw[thick] (#2) -- (#3) -- (#4) -- cycle;  
}
\newcommand{\PP}{\mathbb{P}}
\newcommand{\pp}{\ensuremath{\mathbf{p}}}
\newcommand{\uu}{\ensuremath{\mathbf{u}}}
\newcommand{\RR}{\ensuremath{\mathbb{R}}}
\newcommand{\IRR}{\ensuremath{\mathbb{IR}}}
\newcommand{\FF}{\ensuremath{\mathbb{F}}}
\newcommand{\QQ}{\ensuremath{\mathbb{Q}}}
\newcommand{\ZZ}{\ensuremath{\mathbb{Z}}}
\newcommand{\ww}{\ensuremath{\mathbf{w}}}
\newcommand{\vv}{\ensuremath{\mathbf{v}}}
\newcommand{\rr}{\ensuremath{\mathbf{r}}}
\newcommand{\xx}{\ensuremath{\mathbf{x}}}
\newcommand{\yy}{\ensuremath{\mathbf{y}}}
\newcommand{\zz}{\ensuremath{\mathbf{z}}}

\newcommand{\qq}{\ensuremath{\mathbf{q}}}

\newcommand{\Ptwo}{\ensuremath{\mathbb{P}^2}}

\newcommand{\GD}{\ensuremath{\operatorname{GD}}}
\newcommand{\conv}{\ensuremath{\operatorname{conv}}}

\newcommand{\DEF}[1]{{\color{purple}\emph{#1}}}

\newcommand{\oSequence}{\sigma}
\newcommand{\intersectionPointOracle}{\Omega_{\text{intersection}}}
\newcommand{\linearApproximationOracle}{\Omega_{\text{linear.approximation}}}
\newcommand{\NewtonOralce}{\Omega_{\text{Newton}}}
\newcommand{\GDRejectingOracle}{\Omega_{\text{GD.disjoint}}}
\newcommand{\GDAcceptingOracle}{\Omega_{\text{GD.converge}}}
\newcommand{\oracle}{\Omega}
\newcommand{\nonzeroCertifiedOracle}{\Omega_{\text{nonzero.certified}}}
\newcommand{\krawczykOralce}{\Omega_{\text{Krawczyk.certified}}}
\newcommand{\triangulation}{\tau}

\newcommand{\safetyCoefficient}{\ensuremath{C_{\text{safety}}}}
\newcommand{\intersectionNormBound}{\ensuremath{C_{\text{max.int.norm}}}}
\newcommand{\areaScalingThreshold}{\ensuremath{C_{\text{area.scaling}}}}
\newcommand{\localTriangle}{\ensuremath{\triangle_{\text{local}}}}

\begin{document}
\title{Projective Plane Subdivision Method\\For Initial
Orbit Determination}
\ifarxiv
\author[1]{Ruiqi Huang}
\author[1]{Anton Leykin}
\author[2]{Michela Mancini}
\affil[1]{School of Mathematics, Georgia Institute of Technology, Atlanta, 30332, GA, USA.}
\affil[2]{Guggenheim School of Aerospace Engineering, Georgia Institute of Technology, Atlanta, 30332, GA, USA.}
\date{August 1, 2025}
\setcounter{Maxaffil}{0}
\renewcommand\Affilfont{\itshape\small}
\fi

\maketitle              
\begin{abstract}
Initial Orbit Determination (IOD) is the classical problem of estimating the orbit of a body in space without any presumed information about the orbit. 
The geometric formulation of the ``angles-only'' IOD in three-dimensional space: find a conic curve with a given focal point meeting the given lines of sight (LOS). 

We provide an algebraic reformulation of this problem and confirm that five is the minimal number of lines necessary to have a finite number of solutions in a non-special case, and the number of complex solutions is~66. 

We construct a subdivision method to search for the normal direction to the orbital plane as a point on the real projective plane. The resulting algorithm is fast as it discovers only a handful of the solutions that are real and physically meaningful. 

\end{abstract}
\section{Introduction}

Initial Orbit Determination (IOD) methods aim at estimating the Keplerian orbit of any negligibly light celestial body around a heavy body (e.g., an asteroid orbiting the Sun, or a satellite orbiting the Earth)

The seminal work of Gauss~\cite{gauss_theoria_1809} addressed, in particular, IOD of Ceres, a then new celestial body observed in the solar system. This work and many others that followed use (in contrast to our approach) \emph{time} as one of the coordinates of an observation.

Removing time --- and therefore dynamics --- from observation data gives us a fascinating purely geometric problem. Noticeably, this approach allows us to solve the problem even in scenarios where time measurements are absent or unreliable (when the classic solutions provided by Laplace ~\cite{laplace_memoire_1780} and Gauss, and the more modern ones provided by Escobal ~\cite{Escobal:1976} and Gooding ~\cite{Gooding:1997} cannot be applied). Moreover, the solution is characterized by its independence from light time-of-flight corrections. These corrections can be substantial and are challenging to estimate when the distance between the observer and the body is unknown, as is the case of the IOD problem.
  
\medskip
The geometric formulation of the angles-only IOD problem consists of finding a conic curve in $\RR^3$ with a focal point at the origin that is incident to five given lines, which are assumed to be generic.\footnote{In the precise language of algebraic geometry a configuration of five lines is \DEF{generic} if it avoids a (perhaps known but fixed) hypersurface in the variety that is the fifth power of the Grassmannian of lines in the three-dimensional space. In practical terms, a generic configuration is non-special with respect to the stated problem or, more narrowly, a method of solving the problem: in particular, the method doesn't fail if a configuration is chosen at random.}
Our algorithmic approach is both algebraic-geometric and numerical-analytic.

The ingredients include a subdivision routine that produces an adaptive triangulation of $\Ptwo$, the projective plane of normal directions to the orbital plane, and a fast evaluation routine for an auto-differentiable map on the projective plane that vanishes exactly on the solutions to the problem. What makes our algorithm fast and practical is the use of \DEF{oracles} guiding the refinement of the triangulation. 
Some of the triangles are marked as ``accept'' --- containing a solution --- or ``reject'' --- containing no solution. In theory, one can provide oracles that certify (formally verify) this. In practice, due to higher computational cost, the certified oracles may be used for post-processing the results obtained by reliable yet approximate oracles that use various relaxations.

Certified approaches are based on using interval arithmetic as well as the Krawczyk method~\cite{krawczyk1969newton} (e.g., explained in \cite{burr2019effective}). Heuristic oracles utilize the insight of numerical approaches like Newton's method and gradient descent, as well as some common-sense physical constraints. 

\medskip

Our approach is direct, as opposed to the dual approach of Duff et al.~\cite{duff2022orbit}, to the same problem. In a nutshell, it boils down to finding a handful of approximate \emph{real} solutions to a system of two (analytic) equations in two unknowns. This is in contrast to finding 66 \emph{complex} solutions of seven (polynomial) equations in seven unknowns in the most optimal dual formulation of~\cite{duff2022orbit}.    

The two methods are very different: the workhorse of~\cite{duff2022orbit} is (complex) polynomial homotopy continuation, while the algorithm of this paper relies on a subdivision method in (real) dimension two. 

It shall be mentioned that the problem of solving (polynomial) systems via subdivision methods is addressed by, for instance,~\cite{mourrain2009subdivision} in full rigor. The value of our approach, given its area of practical application, is in fast heuristic numerical approximation. Nevertheless, we address the question of rigorous post-certification of found approximations in~\Cref{section:certification} and computing the algebraic degree of the problem in~\Cref{section:solution-count}.

\section{Geometric outline of IOD problem}
\begin{wrapfigure}[15]{r}{0.45\textwidth}
    \centering
    \includegraphics[width=0.9\linewidth]{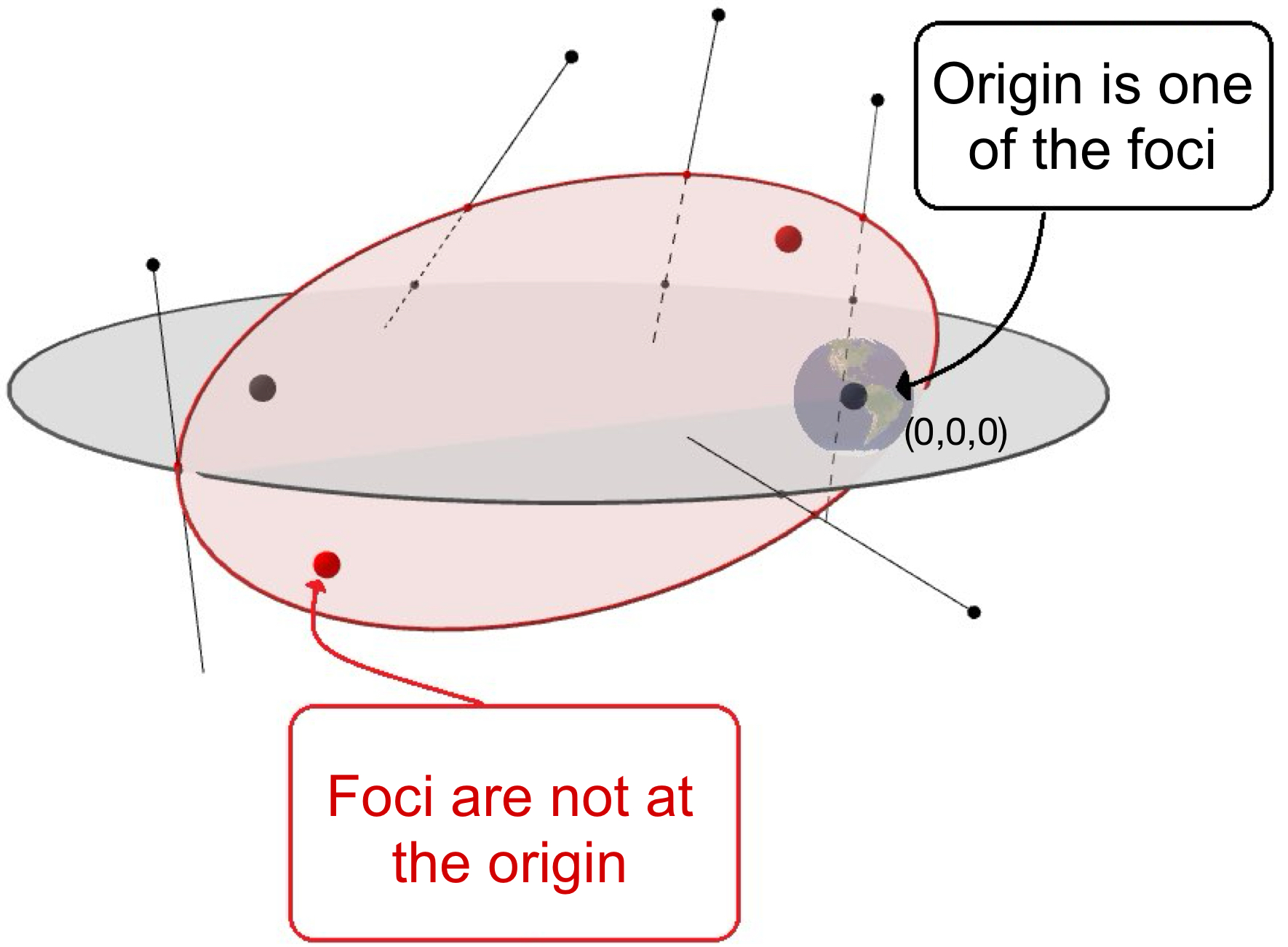}
    \caption{The example and non-example of the orbit satisfying the assumptions}
    \label{fig:concept-illustration}
\end{wrapfigure}

An angle-only observation is composed of the position of the observer and the direction of observation, which define a line in 3-D space that must intersect the orbit of the observed orbiting body. Given five lines in space, and assuming Keplerian dynamics, our objective is to determine the orbit. Under Keplerian dynamics, the orbiting object's orbit is a conic section with a focus at the gravitating body. The location of such
body is assumed to be known; therefore, so is a focal point of the conic. The gray orbit in \Cref{fig:concept-illustration} is an orbit satisfying our assumptions, while the red ellipse is not the right orbit since neither of its focal points is located at the origin.

We note that if the orbital \emph{plane} is known, then finding the conic reduces to a linear algebra problem of finding a quadratic equation that vanishes at five points on the plane, the intersection points of the lines of sight with the known plane. This, in turn, shows that if we know the normal vector to the orbital plane and five lines of sight, we are able to find the conic equation.

We shall assume, without loss of generality, that the heavy body is placed at the origin and think of the normal direction to the plane of the orbit as an indeterminate point on the projective plane $\PP^2$. Thus, given the lines of sight information, typically provided in the form of the position vectors $\pp_i \in \RR^3$ and direction vectors $\uu_i\in \RR^3$ of the observations, for $i=1,\dots,5$, we want to find the vector $\ww \in \PP^2$ normal to the needed orbital plane, assuming the origin is one of its focal points.

\section{Subdivision of the projective plane}

We represent the real projective plane $\Ptwo$ as the surface of the octahedron with vertices at $\pm e_1, \pm e_2, \pm e_3$, where $e_1,e_2,e_3$ is the standard basis of $\RR^3$, and where the antipodal points are identified.

\begin{figure}[H]
  \centering
    \subfloat[\centering]{{\includegraphics[width=5cm]{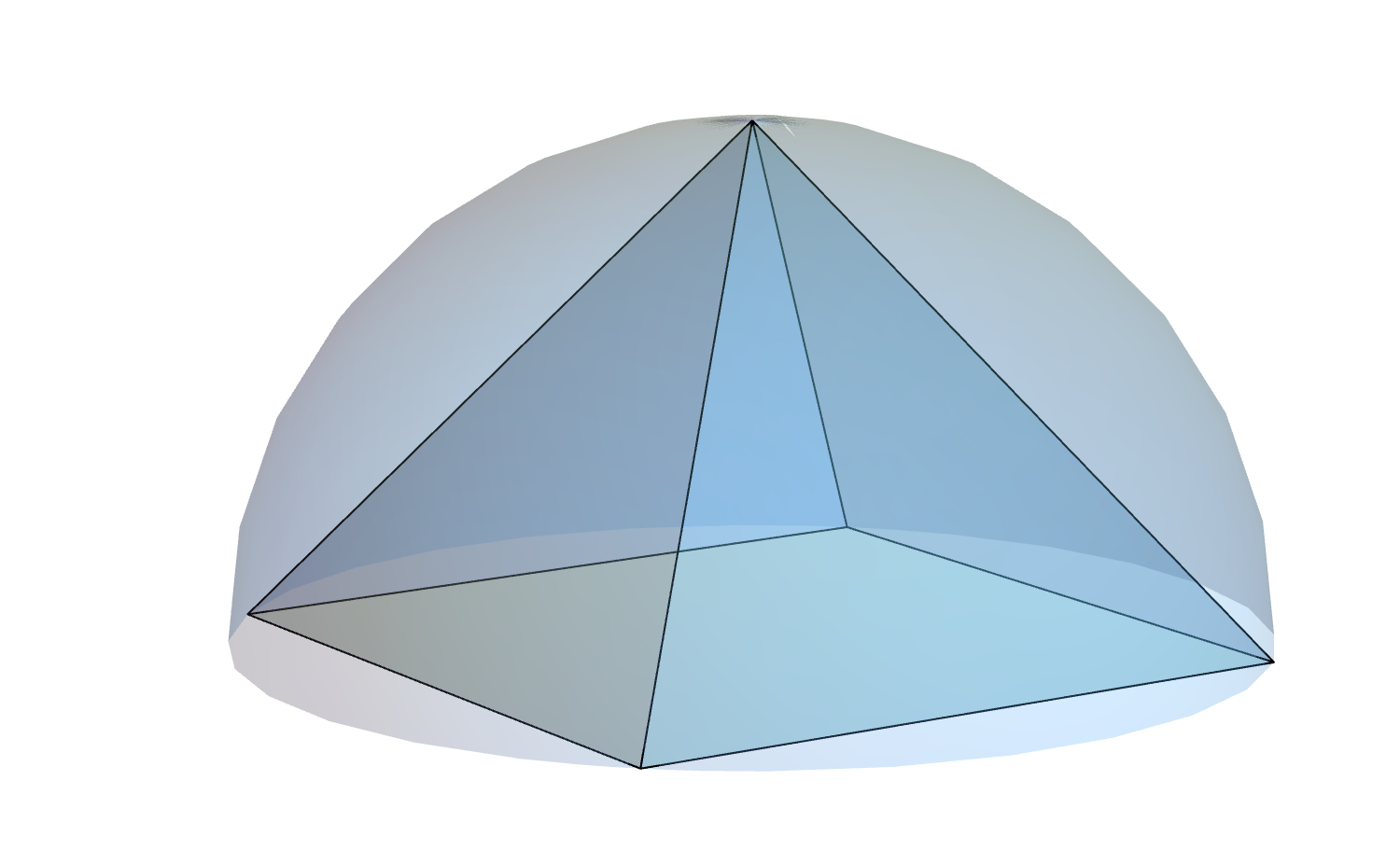} }}%
    \subfloat[\centering]{{\includegraphics[width=3.5cm]{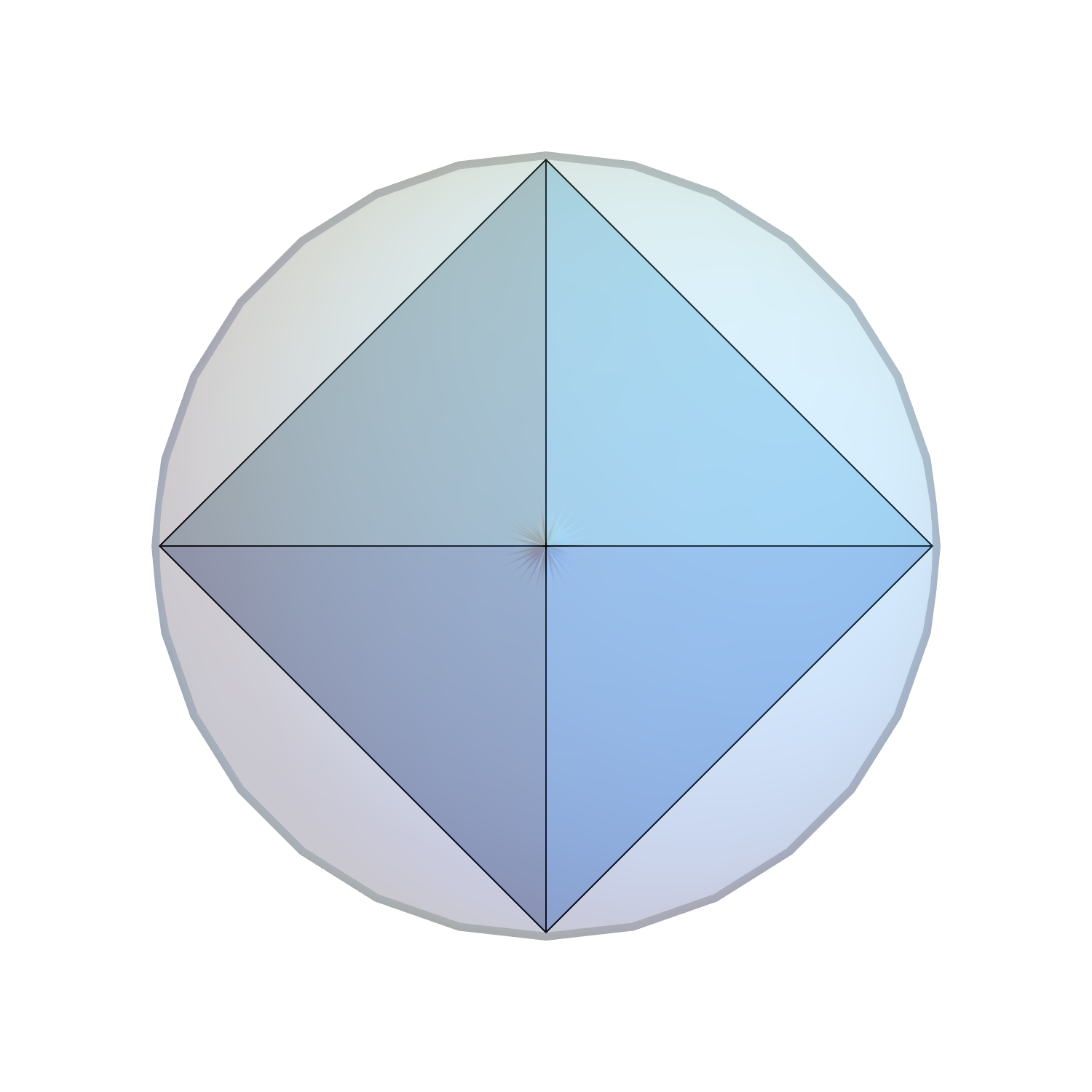} }}%
    \label{figure:half-octahedron}
  \caption{Half of the octahedron representing~$\Ptwo$ viewing from the front (a) and viewing from the above (b).}
\end{figure}

\begin{figure}[H]
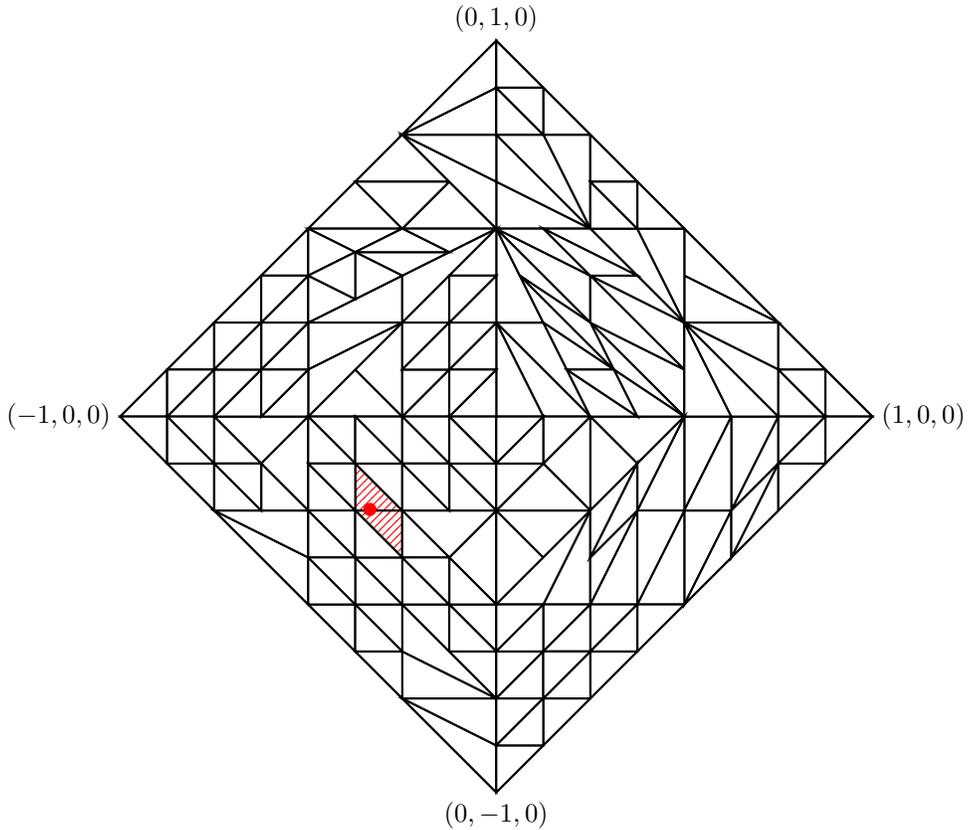

    \centering
    \begin{tikzpicture}[scale=5]
        \input{Tikz_TXT_Files/Skeleton.txt}
        \input{Tikz_TXT_Files/Almost_Parallel.txt}
    \end{tikzpicture}
    \caption{Example of a triangulation obtained for synthetic data with a unique relevant solution (the true solution of $w$ is marked in red, and the triangles shaded in red are triangles accepted by our method)}
    \label{figure:almost-parallel}
\end{figure}

Our algorithm shall subdivide the ``northern'' faces of the octahedron, replacing the four large triangles with unions of several smaller triangles. For example, $\Cref{figure:almost-parallel}$ is an example of the triangulation history after applying the algorithm. In fact, we use subdivision algorithms that restrict the vertices in any triangulation to the set   
\[
\left\{\left(\frac{x}{2^a},\frac{y}{2^b},\frac{z}{2^c}\right) : \left|\frac{x}{2^a}\right|+\left|\frac{y}{2^b}\right|+\left|\frac{z}{2^c}\right| = 1;\, a,b,c,x,y,z\in\ZZ;\, a,b,c,z\geq 0\right\} \subset \QQ^3
\]
thus warranting exactness.
Given a triangle $\triangle$ in the current triangulation, a projective plane subdivision algorithm in this article shall decide whether to 
\begin{itemize}
    \item mark $\triangle$ as \emph{rejected}
    \begin{itemize}
        \item either due to the impossibility of having a solution of the system in $\triangle$,
        \item or due to \emph{reasonable constraints} imposed on the intermediate elements in the construction, for instance, 
        \begin{itemize}
            \item a lower bound on the angle a line of sight forms with the orbital plane,
            \item an upper bound on the norm of an intersection point of a line of sight and the orbital plane,
            \item a bound on sensitivity of the conic fit of the intersection points to noise, etc.;
        \end{itemize}
    \end{itemize}
    \item mark $\triangle$ as \emph{accepted} if (depending on our goal) it is shown that it 
        \begin{itemize}
            \item either contains at least one solution
            \item or contains exactly one solution;
        \end{itemize}
    \item or \emph{subdivide} $\triangle$ into smaller triangles
        \begin{itemize}
            \item either in a predetermined manner (e.g., taking four equal triangles with vertices at the vertices of $\triangle$ or the midpoints of the sides of $\triangle$)
            \item or \emph{adaptively} (e.g., basing the type of subdivision on the local approximation of the Jacobian).
        \end{itemize}
\end{itemize}

Assuming that a solution set is finite (it is infinite only if the problem is degenerate) the subdivision procedure shall terminate with a collection of accepted triangles, and all other triangles in the final triangulation rejected.

\section{Composite function} \label{section:master-function}
\newcommand{\inputSpace}{\mathcal L}
\newcommand{\localCoordinateSpace}{\mathcal X}
\newcommand{\vanishSpace}{\mathcal V}
\newcommand{\normalSpace}{\mathcal N}
\newcommand{\frameSpace}{\mathcal F}
\newcommand{\fivePointsSpace}{\mathcal P}
\newcommand{\conicSpace}{\mathcal C}

To solve our problem, we shall set up a system of two equations in two unknowns, $F=0$, where $F:\RR^2\to\RR^2$ is our ``master'' function.  

For a triangle $\triangle$ and parameters $\ell\in\inputSpace$ specifying the lines of sight, we design the function $F$ as the following composition 
\begin{equation}\label{equation:master-function}
  \begin{tikzcd}
    \vanishSpace &
    \arrow[l, "F_{\vanishSpace\conicSpace}"] 
    \conicSpace &
    \arrow[l, "F_{\conicSpace\fivePointsSpace}"] 
    \fivePointsSpace & 
    \arrow[l, "F_{\fivePointsSpace\frameSpace}"] 
    \frameSpace &
    \arrow[l, "F_{\frameSpace\normalSpace}"] 
    \normalSpace &
    \arrow[l, "F_{\normalSpace\normalSpace}"] 
    \normalSpace &
    \arrow[l, "F_{\normalSpace\localCoordinateSpace}"] 
    \localCoordinateSpace
  \end{tikzcd}
\end{equation}
where
\begin{itemize}[label={\ding{228}}]
\item $\localCoordinateSpace\cong\RR^2$ is the local coordinate space of a triangle $\triangle$,
\item $\normalSpace\cong\RR^3$ is the space where normal (to the orbital plane) vectors are taken,
\item $\frameSpace\cong\RR^9$ is the space for an orthonormal frame (three vectors in $\RR^3$) that completes the normal vector,
\item $\fivePointsSpace\cong\RR^{10}$ is the space of five points of intersection of the lines of sight with the orbital plane (in planar coordinates), 
\item $\conicSpace\cong\RR^5$ is the space of planar conic equations (five coefficients of a quadric; the free term is fixed),
\item $\vanishSpace\cong\RR^2$ is the space of values of the two polynomials in the coefficients of the quadric that vanish if and only if $(0,0)$ is a focus of the conic.  
\end{itemize}

Therefore, the master function $F$ is constructed by taking the function composition:
\begin{equation*}
    F = F_{\vanishSpace\conicSpace}\circ F_{\conicSpace\fivePointsSpace}\circ F_{\fivePointsSpace\frameSpace}\circ F_{\frameSpace\normalSpace} \circ F_{\normalSpace\normalSpace} \circ F_{\normalSpace\localCoordinateSpace}.
\end{equation*}

For our algorithm, we need the evaluation of $F$ and its Jacobian $J$ to be efficient; $J$ is computed by the chain rule, i.e., by evaluating the Jacobians of its parts and then composing them.

\begin{figure}[H]
    \centering
    \begin{tikzpicture}[scale=3]
      \draw[->] (-0.8,0) -- (1.2,0) node[anchor=west] {$\mathrm{}$};
      \draw[->] (0,-0.8) -- (0,0.8) node[anchor=south] {$\mathrm{}$};
    
      \filldraw[black] (1,0) circle (0.5pt) node[anchor=south west] {\footnotesize $\mathrm{Q}\, (1, 0)$};
      \filldraw[black] (-0.5,0.5) circle (0.5pt) node[anchor=south east] {\footnotesize $\mathrm{R}\, (-0.5, 0.5)$};
      \filldraw[black] (-0.5,-0.5) circle (0.5pt) node[anchor=north east] {\footnotesize$\mathrm{P}\, (-0.5, -0.5)$};
      \filldraw[black] (0,0) circle (0.5pt) node[anchor=south] {\footnotesize center};
    
      \draw[thick] (-0.5,0.5) -- (-0.5,-0.5) -- (1,0) -- cycle;
    
    \end{tikzpicture}
    \caption{Any triangle in the local coordinates}
    \label{figure:triangle-in-local-coordinates}
\end{figure}
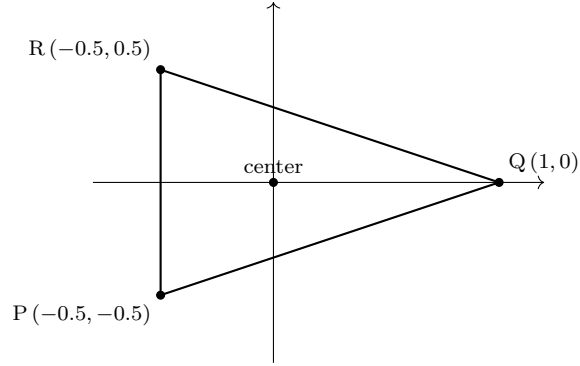

\subsection{Local coordinates}
Given a triangle $\triangle \subset \RR^3$, the first map in the composition provides a parametrization of $\triangle$. It maps the triangle in \Cref{figure:triangle-in-local-coordinates} to $\triangle$ via the linear map $F_{\normalSpace\localCoordinateSpace}$ in \Cref{equation:master-function}. This parametrization makes sure that the center of the triangle in \Cref{figure:triangle-in-local-coordinates} is the origin of the local coordinates.

\subsection{Normalization} \label{section:normalization} 
The second map of~\Cref{equation:master-function} is defined for all points in $\normalSpace \cong \RR^3$ except the origin: 
\[
F_{\normalSpace\normalSpace}(\ww) = \frac{\ww}{\|\ww\|}.
\]

\subsection{Orthonormal frame}
The next map, $F_{\frameSpace\normalSpace}$ takes (a unit vector) $\ww$ and completes it to an orthonormal frame $(\ww,\vv_1,\vv_2)$ where
\begin{equation*}
    \vv_2 = \frac{\ww \times \uu_1}{\lvert \ww \times \uu_1 \rvert}, 
    \quad 
    \vv_1 = \frac{\vv_2 \times \ww}{\lvert \vv_2 \times \ww \rvert}
\end{equation*}
Recall that $\uu_1 \in \RR^3$ is the first of the direction vectors.

\subsection{Intersection with the orbital plane} \label{section:intersectOrbitPlane} 
Given $(\ww, \vv_1, \vv_2)\in\frameSpace$ the map $F_{\fivePointsSpace\frameSpace}$ sends it to $ (\xx,\yy) \in \fivePointsSpace \cong \RR^{10}$, where $(x_i,y_i)$ are $(\vv_1,\vv_2)$-coordinates of the point $\rr_i \in \RR^3 $ where the $i$-th line of sight intersects the plane of the orbit.  

Let $\rr_i$ with $i = 1,\cdots,5$ be the coordinates of the intersection point in $\RR^3$ constructed as $\rr_i = \pp_i + \rho_i\uu_i$ where $\pp_i$ and $\uu_i$ are the $i$-th position and direction vectors of the observer, and $\rho_i$ is the distance from the point to the location of the $i$-th observer. To specify $\rho_i$, we could apply the constraint $\rr_i^T\ww = 0$, which gives
\begin{equation*}
    \rho_i = - \frac{\pp_i^T\ww}{\uu_i^T\ww}.
\end{equation*}

\subsection{Fitting the conic} \label{fitConic} 
Given planar coordinates of five points $(x,y) \in \fivePointsSpace$, the map $F_{\conicSpace\fivePointsSpace}$ gives coefficients $\boldsymbol{\theta} = [a,b,c,d,e]^T \in \conicSpace \cong \RR^5$ in the quadratic equation describing the conic: 
\begin{equation}\label{equation:quadric}
    a x^2 + b xy + c y^2 + d x + e y + 1 = 0.
\end{equation}
Note that for a normal direction that solves our problem, the origin $(0,0)$ has to be a focal point and doesn't satisfy the quadric equation. Thus, specializing the constant term to $1$ is natural. 

Let $M = \left[\begin{matrix}
    \xx^2 & \yy^2 & \xx\yy & \xx & \yy
\end{matrix}\right]$ where
\begin{equation*}
    \begin{aligned}
        \xx^2 &= \left[\begin{matrix}
            x_1^2 & x_2^2 & x_3^2 & x_4^2 & x_5^2
        \end{matrix}\right]^T,\ 
        \yy^2 = \left[\begin{matrix}
            y_1^2 & y_2^2 & y_3^2 & y_4^2 & y_5^2
        \end{matrix}\right]^T, \text{ and }\\
        \xx\yy &= \left[\begin{matrix}
            x_1y_1 & x_2y_2 & x_3y_3 & x_4y_4 & x_5y_5
        \end{matrix}\right]^T.
    \end{aligned}
\end{equation*}
We obtain coefficients $\boldsymbol{\theta}$ by solving the linear system
\begin{equation*}
    M\boldsymbol{\theta} = [\begin{matrix}
        -1 & -1 & -1 & -1 & -1
    \end{matrix}]^T.
\end{equation*}

\subsection{Polynomials that vanish} \label{section:vanishingFreePoints} 
Having $\boldsymbol{\theta} = (a,b,c,d,e) \in\conicSpace$, the coefficients of the quadric in \Cref{equation:quadric}, the final map is \[
F_{\vanishSpace\conicSpace}(\mathbf{\boldsymbol{\theta}}) =  \begin{bmatrix}
        e^2-4b-d^2+4a \\ de-2c
    \end{bmatrix} \in \vanishSpace.
\]
by considering the relation between the focal point and the dual conic~\cite{semple_algebraic_1998}.

\subsection{Alternative composite functions}
We remark that one may want to experiment with other functions
that indicate whether the origin is a focus of the conic.
For instance, one can replace the last two maps in \Cref{equation:master-function}
by
\newcommand{\circularDualConicSpace}{\mathcal D}
\begin{equation}\label{equation:master-function-alt}
  \begin{tikzcd}
    \vanishSpace &
    \arrow[l, "G_{\vanishSpace\circularDualConicSpace}"] 
    \circularDualConicSpace &
    \arrow[l, "G_{\circularDualConicSpace\fivePointsSpace}"] 
    \fivePointsSpace & 
    \arrow[l,""]
    \dots
  \end{tikzcd}
\end{equation}
where
\begin{itemize}[label={\ding{228}}]
\item $\circularDualConicSpace\cong\RR^5$ is the space of values of certain monomials in the coefficients of the dual conic, which has to be a circle in case $(0,0)$ is a focus, and
\item $\vanishSpace\cong\RR^3$ is the space of values of three polynomials that vanish if and only if there is no discrepancy in the monomial values.
\end{itemize}
We have not explored this option experimentally as it is less convenient to set up than the primary approach. 

\section{Upper bound on the number of solutions}
\label{section:solution-count}
\subsection{Equations and unknowns}
With parts of the composition described in the master function (\ref{equation:master-function}) mentioned in \Cref{section:master-function}, one could construct a system of polynomial constraints that ensures the elliptical orbit found has the origin as one of its foci. In total, the construction described below produces a system of 15 polynomial equations in 15 variables.

\subsubsection*{Normalization and Orthonormal frame}
We first introduce the unknowns $w_1,w_2,w_3$, which gives the normal vector $\ww = [w_1 \; w_2 \; w_3]^T$. Then, we could construct the $(\ww,\vv_1,\vv_2)$-frame similarly to \Cref{section:intersectOrbitPlane}. 
Introducing 
\begin{equation*}
    \lambda_2 = \frac{1}{\|\ww \times \uu_1\|}, \quad \lambda_1 = \frac{1}{\|\vv_2 \times \ww\|},
\end{equation*}
we can create unit vectors $\vv_2 = \lambda_2(\ww \times \uu_1)$ and  $\vv_1 = \lambda_1 (\vv_2 \times \ww)$. 

In summary, we introduce $\lambda_1$ and $\lambda_2$ as variables (together with $w_1,w_2,w_3$) constrained by setting squares of the norms of $\vv_1$, $\vv_2$, and $\ww$ equal to $1$.

\subsubsection*{Intersection with the orbital plane}
For $i = 1,\dots,5$, we introduce variables $\rho_i$ that give the intersection points $\rr_i = \pp_i+ \rho_i\uu_i$ of the lines of sight with the orbit plane that satisfy $\rr_i^T\ww = 0$, since the vectors in the orbital plane are orthogonal to the normal vector $\ww$. This adds 5 variables and 5 constraints.

\subsubsection*{Fitting the Conic}
After normalizing the general symmetric matrix $C$ of the quadratic form (by setting the bottom right entry to $1$) we have five new variables $a,b,c,d,e$, which are the coefficients of the equation of an ellipse as described in~\Cref{equation:quadric}.

There are five constraints 
\begin{equation*}
    \rr_i'^TC\rr_i' = 0 \text{ for } i = 1, \dots,5
\end{equation*}
where $\rr_i'$ are computed from $\rr_i$ using the $(\vv_1,\vv_2)$-frame for the orbit plane. 

\subsubsection*{Polynomials that vanish}
Two equations are produced by constraining a focal point to the origin without introducing new variables.

\subsection{Solution count and symmetries}

One could use the 15 constraints above as generators to construct an ideal
\begin{equation*}
    I \subseteq \FF_p[w_1,w_2,w_3,\lambda_1,\lambda_2,\rho_1,\rho_2,\rho_3,\rho_4,\rho_5,a,b,c,d,e].
\end{equation*}
Gr\"{o}bner basis computations modulo several large primes $p$ (i.e. over a finite field $\FF_p$), reveal that this is a zero-dimensional ideal with degree $528$ for a generic choice of input (for many random choices). In geometric terms, this means that $528$ is the number of complex solutions of the system of equations formed by the constraints, which gives an upper bound for the number of real solutions to the system.

However, since some symmetries are created by the constraints, the upper bound for the solutions to the IOD problem could be less than $528$. Notice for the constraint $\lVert \ww \rVert^2=1$, if $(w_1,w_2,w_3)$ is a solution, then $(-w_1,-w_2,-w_3)$ could also be a solution with the other variables unchanged. Similarly, considering the constraints $\lVert \vv_i \rVert^2 =1$ for $i = 1,2$, one can see that if $\lambda_i$ is a solution, so is $-\lambda_i$. Then, for any solution, we have $2\times 2 \times 2 = 8$ solutions created by the symmetries representing the same elliptical orbit. Thus, we should have $\tfrac{528}{8} = 66$ as an upper bound of the real solutions to the angle-only IOD problem. This confirms the count produced with the dual method in~\cite{duff2022orbit}.

\section{Oracles: reject, accept, or subdivide?}

Here we shall describe the several methods, which we refer to as \DEF{oracles}, that, given a triangle $\triangle$, make a decision whether to 
\begin{itemize}[label={\ding{228}}]
    \item \DEF{``reject''},
    \item \DEF{``accept''}, or
    \item \DEF{``pass''}, that is, make no conclusion. 
\end{itemize}
Our implementation of the subdivision method depends on the choice of an \DEF{oracle sequence} of length $N_\oracle$, 
\[ 
\oSequence = \left(\oracle_1, \dots, \oracle_{N_\oracle}\right),
\]
and for every $\triangle$ in the current triangulation $\triangulation$ applies a simple \Cref{alg:label-triangle}.
\begin{algorithm}[htbp]
\caption{LabelTriangle}\label{alg:label-triangle}
\begin{algorithmic}
    \Procedure{LabelTriangle}{$\triangle,\oSequence$}
    \State $L \gets \text{``pass''}$
    \State $i \gets 1$ 
    \While{$L=\text{``pass''} \And i\leq N_\oracle$}
        \State $L \gets \oracle_i(\triangle)$
        \State $i \gets i+1$
    \EndWhile
    \State\Return $L$
    \EndProcedure
\end{algorithmic} 
\end{algorithm}

This algorithm is designed by the idea that every triangle will be labeled as ``pass'' at first, and then based on each oracle $\oracle_i$ in the oracle sequence $\oSequence$, the label of the triangle might be changed to ``accept'' or ``reject'', and if the oracles are not able to change the label, it will remain as ``pass''.

Below are descriptions of oracles based on various principles or methods.

\subsection{Intersection points feasibility ($\intersectionPointOracle$)}

This oracle $\intersectionPointOracle$ labels $\triangle$ as ``reject'' if the computation of points of intersection between the lines of sight and the orbital plane (\Cref{section:intersectOrbitPlane}) performed at the center of this triangle delivers points on the orbital plane with $\ell_2$-norm exceeding $\intersectionNormBound$. 

\subsection{Linear approximation ($\linearApproximationOracle$)}

Starting with the local parametrization of a triangle $\triangle$ as in \Cref{figure:triangle-in-local-coordinates}, the linear approximation oracle $\linearApproximationOracle$ applies the ``reject'' label to a $\triangle$ if at the center $\mathbf{0} = (0,0)$ of the triangle in the local coordinates.
\begin{equation*}
    \lVert F\left(\mathbf{0}\right) \rVert_2 - \safetyCoefficient\lVert J_{\mathbf{0}} \rVert > 0
\end{equation*}
 where $\safetyCoefficient$ is a constant set by the user, and $\lVert J_{\mathbf{0}} \rVert$ is the matrix norm of the Jacobian at the center of the triangle. This oracle is derived from the first-order Taylor expansion centered at zero with the perturbation size of the safety coefficient $\safetyCoefficient$. This means if the norm of the master function evaluated at all of the points that are within the circle centered at the center of the triangle with radius $\safetyCoefficient$ fails to obtain zero based on the linear approximation by the matrix norm of the Jacobian, then it is less likely to have a true solution inside the triangle. Since the circle with radius $\safetyCoefficient$ covers the triangle, and there is no zero inside this ball by the approximation, the oracle decides to ``reject''.

\subsection{Newton's method ($\NewtonOralce$)}
Given a triangle $\triangle$, we denote the triangle in \Cref{figure:triangle-in-local-coordinates} by $\localTriangle$. Label the vertices of $\localTriangle$ as $P,Q$ and $R$. For a function $f$, we define the Newton operator applied to a vector $\zz$ to be
\begin{equation*}
    N_f(\zz) = \zz - J_\zz^{-1}f(\zz).
\end{equation*}
Let $P' = N_F(P)$, where $F$ is the master function, and similarly let $Q'$ and $R'$ be the points after applying the Newton operator to $Q$ and $R$. Then, denote the triangle formed by $P',Q'$ and $R'$ in the local coordinates as $\localTriangle'$. The Newton oracle $\NewtonOralce$ applies the ``accept'' label to $\triangle$ if
\begin{itemize}
    \item $\localTriangle' \subset \localTriangle$ and
    \item Area$(\localTriangle') \leq \areaScalingThreshold \cdot$ Area$(\localTriangle)$
\end{itemize}
where $\areaScalingThreshold$ is a parameter chosen by the user. This means Newton's method moves the vertices of the triangle inward and refines the region possibly containing solutions to a comparatively smaller triangle, which indicates the possibility of the existence of the zero of the master function inside the triangle.

However, one can always modify the procedure by sampling a different set of points on $\localTriangle$. For example, one can apply the Newton operator to not only the three vertices but also the midpoints of the sides $m_{PQ}, m_{QR}$, and $m_{PR}$. Then, assuming the points after applying the Newton operator to them are $m_{PQ}',m_{QR}'$ and $m_{PR}'$, one could create the polygon $H = \conv(P',Q',R',m_{PQ}',m_{QR}',m_{PR}')$ by finding the convex hull of the six points to make a more reliable decision than using only three points. $\NewtonOralce$ labels a $\triangle$ by ``accept'' if
\begin{itemize}
    \item $H \subset \localTriangle$ and
    \item Area$(H)$$ \leq \areaScalingThreshold \cdot$ Area$(\localTriangle)$.
\end{itemize}

\subsection{Gradient descent method ($\GDAcceptingOracle$ and $\GDRejectingOracle$)}
\begin{wrapfigure}{R}{0.4\textwidth}
    \centering
    \begin{tikzpicture}[scale=2.5]
    
      \coordinate (P) at (-0.5,-0.5);
      \coordinate (Q) at (1,0);
      \coordinate (R) at (-0.5,0.5);
      \draw[thick] (P) -- (Q) -- (R) -- cycle;
    
      \coordinate (O) at (0,0);
      \filldraw[black] (O) circle (0.4pt) node[below left] {$O$};
    
      \coordinate (G) at (0.6,0.5);
      \filldraw[black] (G) circle (0.2pt);
    
      \path[name path=gradient] (O) -- (G);
      \path[name path=edgeQR] (Q) -- (R);
    
      \path[name intersections={of=gradient and edgeQR, by=T}];
      \fill (T) circle (0.4pt) node[above] {$T$};
    
      \draw[->, thick] (O) -- (G) node[above] {$\nabla f$};
    
      \coordinate (M) at ($0.5*(O)+0.5*(T)$);
      \fill (M) circle (0.4pt) node[below right] {$M$};
    
    \end{tikzpicture}
    \caption{Points used to choose the learning rate}
    \label{fig:gradient}
\end{wrapfigure}
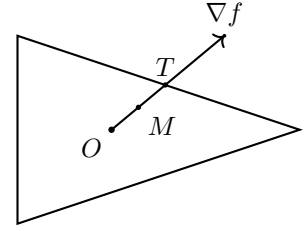
For a function $f$, we define the one-step gradient descent operator $GD_f$ to be
\begin{equation*}
    \GD_f(\zz) = \zz - \gamma \nabla f (\zz)
\end{equation*}
where $\gamma$ is the learning rate. Tracing from the center $O$ in the gradient direction, let $T$ be the intersection with the boundary of the triangle. Denote the midpoint of the line segment connecting $O$ and $T$ as $M$, as shown in the \Cref{fig:gradient}. 

We define the learning rate as 
\begin{equation*}
    \gamma = \frac{\lvert(M-O)^T[\nabla f(M)-\nabla f(O)]\rvert}{\lVert \nabla f(M)-\nabla f(O) \rVert_2^2}
\end{equation*}
derived using the Barzilai-Borwein method~\cite{barzilai_two-point_1988}. For any $\qq \in \localCoordinateSpace$, since $F(\qq) = 0$ only if $\lVert F(\qq) \rVert_2^2 = 0$, define $g(\qq) = \lVert F(\qq) \rVert_2^2$ and therefore $\nabla g = 2J_{\qq}F(\qq)$. 

Similarly as in the Newton oracle, starting with a $\triangle$ in the triangulation, we label the vertices of the triangle $\localTriangle$ in the local coordinates as $P$, $Q$, and $R$, and then we apply the gradient descent operator $\GD_g$ to the vertices $P$, $Q$, and $R$ to obtain $P' = \GD_g(P)$, $Q' = \GD_g(Q)$, and $R' = \GD_g(R)$. Then, by considering the triangle $\localTriangle'$ with vertices $P',Q'$ and $R'$, the accepting oracle $\GDAcceptingOracle$ is designed to accept this $\triangle$ if 
\begin{itemize}[label={\ding{228}}]
    \item $\localTriangle' \subset \localTriangle$ and
    \item Area$(\localTriangle') \leq \areaScalingThreshold \cdot$ Area$(\localTriangle)$,
\end{itemize}
which means $\GD_g$ maps the three vertices inside $\localTriangle$ justifying the heuristic acceptance of this triangle.

In contrast, we also design a rejecting oracle $\GDRejectingOracle$ using the gradient descent method, which labels $\triangle$ as ``reject'' if 
\begin{equation*}
    \localTriangle' \cap \localTriangle = \emptyset.
\end{equation*}
This means, heuristically, that $\localTriangle$ is unlikely to contain any solution.

\subsection{Certified methods ($\nonzeroCertifiedOracle$ and $\krawczykOralce$)} \label{section:certification}

To prove the existence and uniqueness of a solution, we define two oracles based on the interval arithmetic and the Krawczyk method.

Given a square system $F: \RR^n \to \RR^n$, we define an interval enclosure $\Box F: \IRR^n \to \IRR^n$ which extends the computations with the interval inputs such that $\{F(a) \lvert a \in A\} \in \Box F(A)$ for all $A \in \IRR^n$. Similarly, $\Box J: \IRR^n \to \IRR^{n\times n}$ is the interval enclosure for the Jacobian of $F$. Let $I \subseteq \IRR^n$ be an interval set, $x_0 \in I$ be a point in $I$, and $[x_0] \subseteq I$ be an interval constructed from $x_0$ with zero width. We define the Krawczyk operator~\cite{krawczyk1969newton} as follows.
\begin{equation*}
    K_{x_0,Y}(I) := x_0 - Y \cdot \Box F([x_0])+(\mathbf{1}_n - Y\cdot \Box J(I))(I-x_0),
\end{equation*}
where $\mathbf{1}_n$ is the $n\times n$ identity matrix and $Y$ is an invertible matrix. Now, we design an oracle based on the following result (stated, for instance, as Theorem 2.1 in~\cite{duff_certified_2024}).

\begin{theorem}\label{theorem:Krawczyk}
    Suppose that $\Box F: \IRR^n \to \IRR^n$ is an interval enclosure based on a square differentiable system. $\Box J: \IRR^n \to \IRR^{n\times n}$ is the interval extension of the Jacobian $J$ of $F$. Let $I \subseteq \IRR^n$, then for an $n \times n$ invertible matrix $Y$ and a point $x_0 \in I$,
    \begin{enumerate}
        \item (Existence) if $K_{x_0,Y}(I) \subset I$, then $I$ contains a solution $x^{\star}$ of $F$,
        \item (Uniqueness) if $\lVert \mathbf{1}_n - Y\cdot \Box J(I) \rVert < 1$, then the solution $x^{\star}$ in $I$ is unique.
    \end{enumerate}
\end{theorem}

Note that the interval matrix norm of an interval matrix $M \in \IRR^{n \times n}$ is defined by
\begin{equation*}
    \lVert M \rVert := \max_{A \in M} \max_{x \in \RR^n} \frac{\lVert Ax \rVert_{\infty}}{\lVert x \rVert_{\infty}}.
\end{equation*}

For our problem, we extend the square system of equations defined by the master function $F: \RR^2 \to \RR^2$ and its Jacobian $J$ to take interval inputs, which are $\Box F: \IRR^2 \to \IRR^2$ and $\Box J: \IRR^2 \to \IRR^{2\times 2}$. Now, for each triangle $\triangle$, we fit the parametrized triangle in the local coordinates into a small box $I \subseteq \IRR^2$ such that it is covered entirely by $I$. Then, we design a rejecting oracle $\nonzeroCertifiedOracle$ that labels $\triangle$ as ``reject'' if 
\begin{equation*}
    0 \notin \Box F(I)
\end{equation*}
Next, we pick $x_0$ to be a point in $I$ (i.e., the midpoints of the intervals) and fix an invertible matrix $Y$. Then, our certified accepting oracle $\krawczykOralce$ accepts $\triangle$ if 
\begin{equation*}
    K_{x_0,Y}(I) \subseteq I.
\end{equation*}
Moreover, if a rigorous uniqueness of a solution is desired, one can use the second part of~\Cref{theorem:Krawczyk}.

\begin{remark} There are a few propositions in the vein of~\Cref{theorem:Krawczyk} involving Krawczyk and Hansen-Sengupta operators in~\cite{neumaier_interval_1991}: for instance, compare Theorem~5.1.8 to 5.1.10 there. The former relied on the computation of the so-called \emph{Lipschitz matrix}. We believe developing a method to estimate this matrix would be harder than computing the spectral radius of the interval matrix $\mathbf{1}_n - Y\cdot \Box J(I)$ which is sufficient for uniqueness in our formulation. It would also likely lead to a more conservative approach.
\end{remark}

\section{Subdivision Methods}
After deciding the labels for all of the triangles in a triangulation $\triangulation$, we stop subdividing or modifying the triangles with labels ``accept'' or ``reject''. 

Then, we subdivide the remaining triangles in $\triangulation$ with label ``pass'' to create a new triangulation and repeat the process of labeling the triangles in this new triangulation. We used a combination of the \textit{regular subdivision} method and the \textit{adaptive subdivision} method to subdivide the triangles aligning with the directions of changes in the master function value.
\subsection{Regular subdivision}

The \textit{regular subdivision} method is a predetermined subdivision method of a triangle. This method first finds the midpoint of the three sides of the triangle and then connects the three midpoints to create four smaller triangles.

\subsection{Adaptive subdivision}
\label{sec:AdaptiveSubdivision}

\begin{wrapfigure}[13]{R}{0.55\linewidth}
    \centering
    \begin{tikzpicture}[scale=3.5, font=\small]
        \coordinate (A) at (0, 0);
        \coordinate (B) at (1, 0);
        \coordinate (C) at (0.5, 0.866); 
        \filldraw (A) circle (0.3pt) node[left = 2pt] {$V_2$};
        \filldraw (B) circle (0.3pt) node[right=2pt] {$V_3$};
        \filldraw (C) circle (0.3pt) node[right=2pt] {$V_1$};
        \coordinate (m1) at ($(A)!0.5!(B)$);
        \coordinate (m2) at ($(B)!0.5!(C)$);
        \coordinate (m3) at ($(C)!0.5!(A)$);
        \draw[thick] (A) -- (B) -- (C) -- cycle;
        \draw[->, thick] (m1) ++(0,0) -- ++(0.2,0) node[below right=1pt] {$\mathbf{t}_1$};
        \draw[->, thick] (m2) ++(0,0) -- ++(-0.100002,0.173204) node[right=1pt] {$\mathbf{t}_2$};
        \draw[->, thick] (m3) ++(0,0) -- ++(-0.100002,-0.173204) node[left =1pt] {$\mathbf{t}_3$};
        \filldraw (m1) circle (0.3pt) node[above = 2pt] {$m_1$};
        \filldraw (m2) circle (0.3pt) node[left=2pt] {$m_2$};
        \filldraw (m3) circle (0.3pt) node[right=2pt] {$m_3$};
        \draw[dashed] (A) -- ($(A)+(0,-0.1)$);
        \draw[dashed] (m1) -- ($(m1)+(0,-0.1)$);
        \draw[<->] ($(A)+(0,-0.1)$) -- node[below] {$d_{\text{side}}/2$} ($(m1)+(0,-0.1)$);
    \end{tikzpicture}
    \caption{Triangle $\triangle$ to be subdivided with the adaptive subdivision method}
    \label{fig:adaptiveSubdivisionIllustration}
\end{wrapfigure}
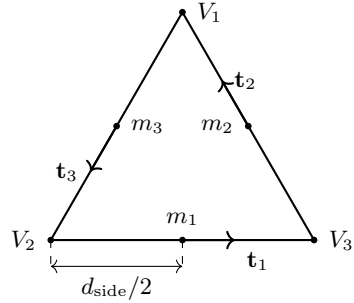

In contrast to the regular subdivision method, the \textit{adaptive subdivision} is a method that takes the Jacobian and master function into consideration to better incorporate with the oracle methods.

As illustrated in the \Cref{fig:adaptiveSubdivisionIllustration}, assume we have a triangle $\triangle$, we also start with finding the midpoints $m_1$, $m_2$, and $m_3$. Then, for the $i$-th side we also compute the unit directional vectors $\mathbf{t}_i$ pointing towards either vertex of this side, with $i = 1,2,3$. Then, we compute the following metric $\delta_i$ for $i= 1,2,3$, to measure the amount of changes in the master function by the linear approximation using the Jacobian $J_{m_i}$ at the midpoints $m_i$.
\begin{equation*}
    \delta_i = \frac{d_\text{side}}{2}\left\lVert J_{m_i} \mathbf{t}_i \right\rVert_2
\end{equation*}
If we have a larger $\delta$ along one side of the triangle, then there will be more variation in the master function value along this side. Therefore, we would connect the midpoint of this side with the opposite vertex to subdivide the triangle into two smaller triangles with the same area. By subdividing along this side, we will better distinguish the triangles with different values of the master functions, which will be helpful to reject or accept the triangular regions with oracles like the linear approximation oracle.

\subsection{Combination of regular and adaptive subdivision methods}

For each triangle labeled as ``pass'' in the triangulation, we first measure the variation in the master function value along the three sides by computing the metric $\delta$ we defined above. Denote the minimum and maximum among $\delta_i$ for $i=1,2,3$ for a triangle as $\delta_{\min}$ and $\delta_{\max}$. Then, if $\delta_{\max} \geq \gamma\cdot \delta_{\min}$, then we apply the adaptive subdivision along the side with $\delta_{\max}$ to make it into two smaller triangles. Otherwise, we apply the regular subdivision into four triangles. The constant $\gamma > 2$ is chosen by the user. We set $\gamma = 4$ for the experiments in the next section.

\section{Data and experiments}

We provide two proof-of-concept examples of the algorithm, and then apply the algorithm to the solution of the IOD problem of a nearly circular orbit with Earth-based observers.

\subsection{Generating proof-of-concept examples}

To generate synthetic data consisting of a normal vector $\ww \in \RR^3$, a position matrix $\pp \in \RR^{3 \times 5}$, and a direction matrix $\uu \in \RR^{3\times 5}$, we could start with fixing the length of the major semi-axis $a$ and the eccentricity $e$ of an ellipse in the $xy$-plane. Then we could recover the distance $c$ between the center and the focal points and the length of the minor semi-axis $b$ by 
\begin{equation*}
    c = a\cdot e, \quad b = \sqrt{a^2-c^2}.
\end{equation*}
Now, we could sample five points $\rr_i \in \RR^3$ by choosing five angles $\beta_i$ measured counterclockwise from the positive $x$-axis, with the modified standard parametrization of the points on an ellipse:
\begin{equation*}
    \rr_i = (a\cos(\beta_i)-c,b\sin(\beta_i),0)
\end{equation*}
since we require one of the focal points to be at the origin. These five vectors give the columns of the matrix $\rr$. For this ellipse, the normal vector is $[\begin{matrix}
    0 & 0 & 1
\end{matrix}]^T$, and once we choose five observer positions to create matrix $\pp^\star$, we can find the direction matrix by $\uu^\star = \rr-\pp^\star$. Finally, by fabricating a rotation matrix $R \in \mathrm{SO}(3)$, we could create a set of fabricated data:
\begin{equation*}
    \ww = R\left[\begin{matrix}
        0& 0 &1
    \end{matrix}\right]^T, \quad \pp = R\cdot \pp^\star, \text{ and } \uu = R \cdot \uu^\star
\end{equation*}

Note that to ensure that we have a normal vector on the ``northern'' faces of the octahedron, we take $\ww = R\cdot [\begin{matrix}
    0 & 0 & -1
\end{matrix}]^T$ if $R \cdot [\begin{matrix}
    0 & 0 & 1
\end{matrix}]^T$ has a negative third coordinate.

If we want to generate synthetic data with two known solutions $\ww_1$ and $\ww_2$ with the same parameters. we could fix two lengths $a_1$ and $a_2$ of major semi-axes with two corresponding eccentricities $e_1$ and $e_2$ and sample points on the ellipses to form matrix $\rr$ and $\rr'$. Then, we first define the rotation matrix $R_2 \in \mathrm{SO}(3)$ to rotate the points $\rr'$ and the true solution $[0 \; 0 \; 1]^T$ from the first one. Then, define the direction matrix $\uu^{\star} = \rr - R_2 \cdot \rr'$ and sample the position vectors in $\pp^{\star}$ as points on the lines defined by the directions. Finally, we could apply the rotation matrix $R_1 \in \mathrm{SO}(3)$ to both solutions and $\pp^{\star}$ and $\uu^{\star}$ to give the pair of solutions and parameters.

\subsection{Setting constants for heuristic oracles}

By the design of the oracles, the constants $\intersectionNormBound$, $\areaScalingThreshold$, and $\safetyCoefficient$ could be chosen by users to make the algorithm do the job that users want. For example, $\intersectionNormBound$ should be set as the distance from the observer to the furthest observable object based on the real physical setting of the problem. On the other hand, the constants should also be decided based on how conservative or aggressive the users want the oracles to be. For example, normally one should have $0 < \areaScalingThreshold < 1$, and $1 \leq \safetyCoefficient$, but if $\areaScalingThreshold$ is $0.5$ or smaller, the acceptance by the Newton oracle would be more conservative, and if $\safetyCoefficient$ is $2$ or higher, the rejection by the linear approximation oracle would be more conservative.

In practice, given a family of problems for a particular physical scenario, one should calibrate the constants in the method to ensure that a true solution is discovered for a few problems sampled from the family.

\subsection{Proof-of-concept examples}

\subsubsection*{Single observer fabricated data}

In this scenario, we choose $\pp$ to be a matrix with its five columns being the same. The synthetic data generated is
\begin{equation*}
    \ww = \left[\!\begin{matrix}
     -{.18511}\\
     -{.944226}\\
     {.272346}
     \end{matrix}\right]
\end{equation*}
\begin{equation*}
         \pp = \left[\begin{matrix}
     {.190367}&{.190367}&{.190367}&{.190367}&{.190367}\\
     -{1.19796}&-{1.19796}&-{1.19796}&-{1.19796}&-{1.19796}\\
     -{.352143}&-{.352143}&-{.352143}&-{.352143}&-{.352143}
     \end{matrix}\right]
\end{equation*}
\begin{equation*}
         \uu = \left[\begin{matrix}
     -{.247813}&{.635429}&{.456836}&-{.536783}&-{.972279}\\
     {1.07074}&{1.07438}&{1.42462}&{1.63744}&{1.41873}\\
     -{.127977}&{.484975}&{1.57787}&{1.64036}&{.58609}
     \end{matrix}\right]
\end{equation*}

With this fabricated data, we applied the oracle sequence
\begin{equation}\label{eqn:OracleScheme}
    \oSequence = \{\intersectionPointOracle,\linearApproximationOracle,\GDRejectingOracle,\NewtonOralce\}
\end{equation}
with the tuned constants for the oracles
\begin{equation*}
    \intersectionNormBound = 10, \quad \areaScalingThreshold = 0.9, \quad \safetyCoefficient = 1.0.
\end{equation*}

As a result, we are able to find one accepted triangle as shown in \Cref{figure:single-observer-result} where the accepted triangles are shaded in red with the true solution (red point) in it, rejected triangles are in white, and blue triangles are triangles with the ``pass'' label.

Practically, the user will set the thresholds for the area of the triangles to activate and terminate the search. For example, in this experiment, the maximal area bound to start applying the oracles is $0.05$, and the minimal area bound to stop the subdivision is $10^{-3}$. This means we will keep subdividing without applying the oracles to label any triangle until the area of a triangle is smaller than the maximal area bound, and we will stop subdividing the triangles once we have the areas of all triangles smaller than the minimal area bound. The result of this experiment shows that the known solution (the red point) is in the accepted triangle (shaded in red) found.

\ifarxiv
\begin{figure}[H]
    \centering
    \includegraphics[width=0.5\linewidth]{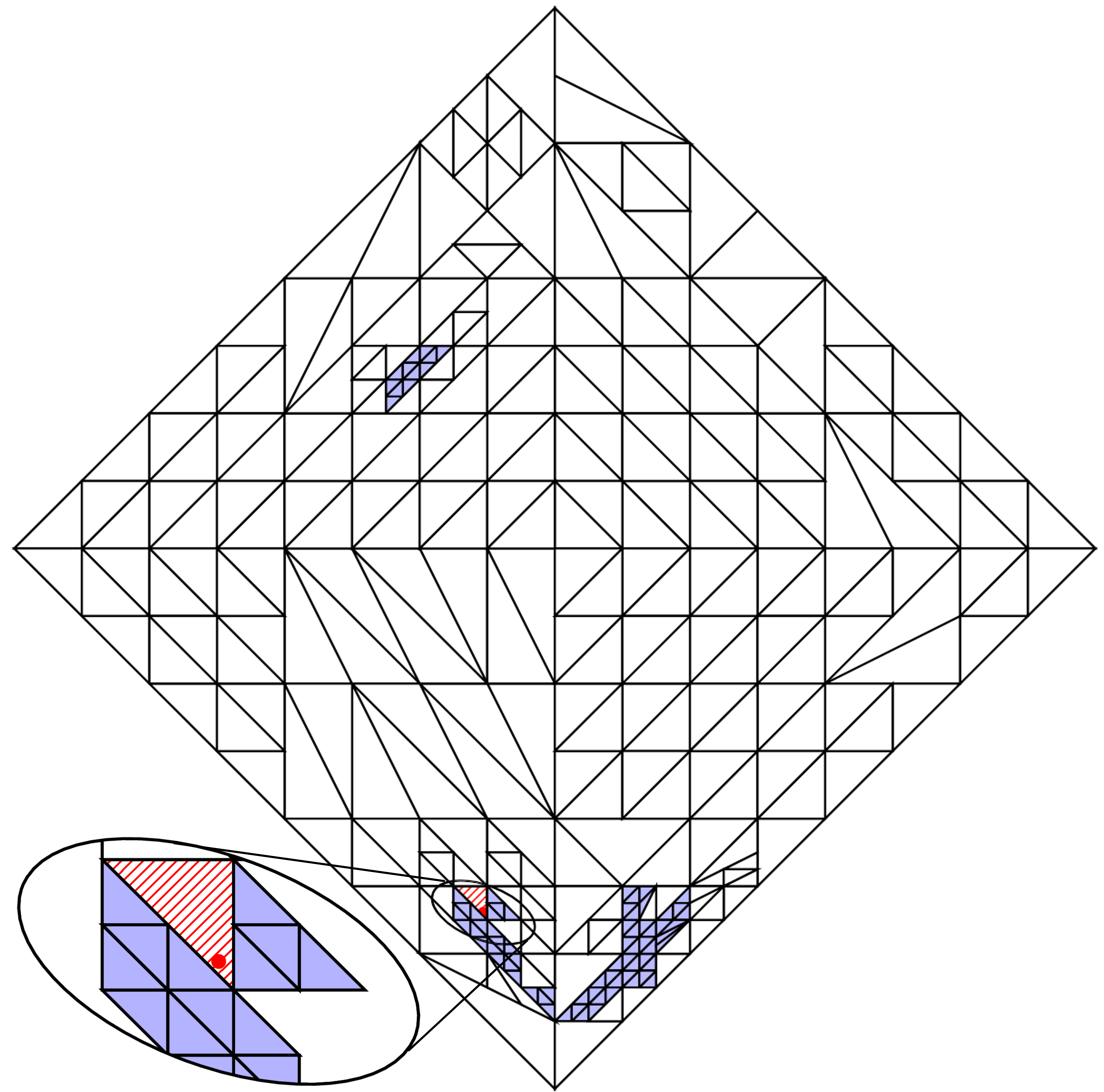}
    \caption{The triangulation from the single observer data with the zoomed-in elliptical region near the known solution}
    \label{figure:single-observer-result}
\end{figure}
\fi

From \Cref{figure:single-observer-result}, we could see that there are also two regions with excessive subdivision (clusters of triangles in blue), which indicates we could possibly accept triangles in these regions if we relax the accepting oracles or subdivide further, and the triangle including the true solution at the bottom-left is zoomed in for better visualization. The zoomed-in region is an ellipse centered at the true solution with the largest and smallest singular values of the Jacobian matrix at the center being the major and minor axes.

\subsubsection*{Two solutions fabricated data}
With the process described above, we could construct a scenario with two fabricated solutions $\ww_1$ and $\ww_2$:
\begin{equation*}
    \ww_1 = \left[\begin{matrix}
      -{.628302}\\
      -{.311317}\\
      {.712964}
      \end{matrix}\right], \quad \ww_2 = \left[\begin{matrix}
      -{.576837}\\
      {.0266409}\\
      {.816425}
      \end{matrix}\right]
\end{equation*}
and the corresponding constants $\pp$ and $\uu$
\begin{equation*}
     \pp = \left[\begin{matrix}
      {.252758}&-{.565296}&-{1.50675}&-{1.27055}&-{.183111}\\
      -{.209549}&-{.906674}&-{2.17693}&-{2.26487}&-{1.04896}\\
      -{.460245}&-{.123196}&{.696834}&{.866591}&{.151477}
      \end{matrix}\right],
\end{equation*}
\begin{equation*}
    \uu = \left[\!\begin{matrix}
      -{.134482}&{.605848}&{.972423}&{.45865}&-{.225455}\\
      -{.121519}&{.420599}&{1.60726}&{1.79853}&{.730092}\\
      {.124171}&{.332124}&{.0711503}&-{.298093}&-{.265324}
      \end{matrix}\right].
\end{equation*}

We applied the same oracle sequence $\oSequence$ described in \Cref{eqn:OracleScheme}, and the following constants for the oracles:
\begin{equation*}
    \intersectionNormBound = 10, \quad \areaScalingThreshold = 0.8, \quad \safetyCoefficient = 0.7.
\end{equation*}
\ifarxiv
\begin{figure}[H]
    \centering
    \includegraphics[width=0.5\linewidth]{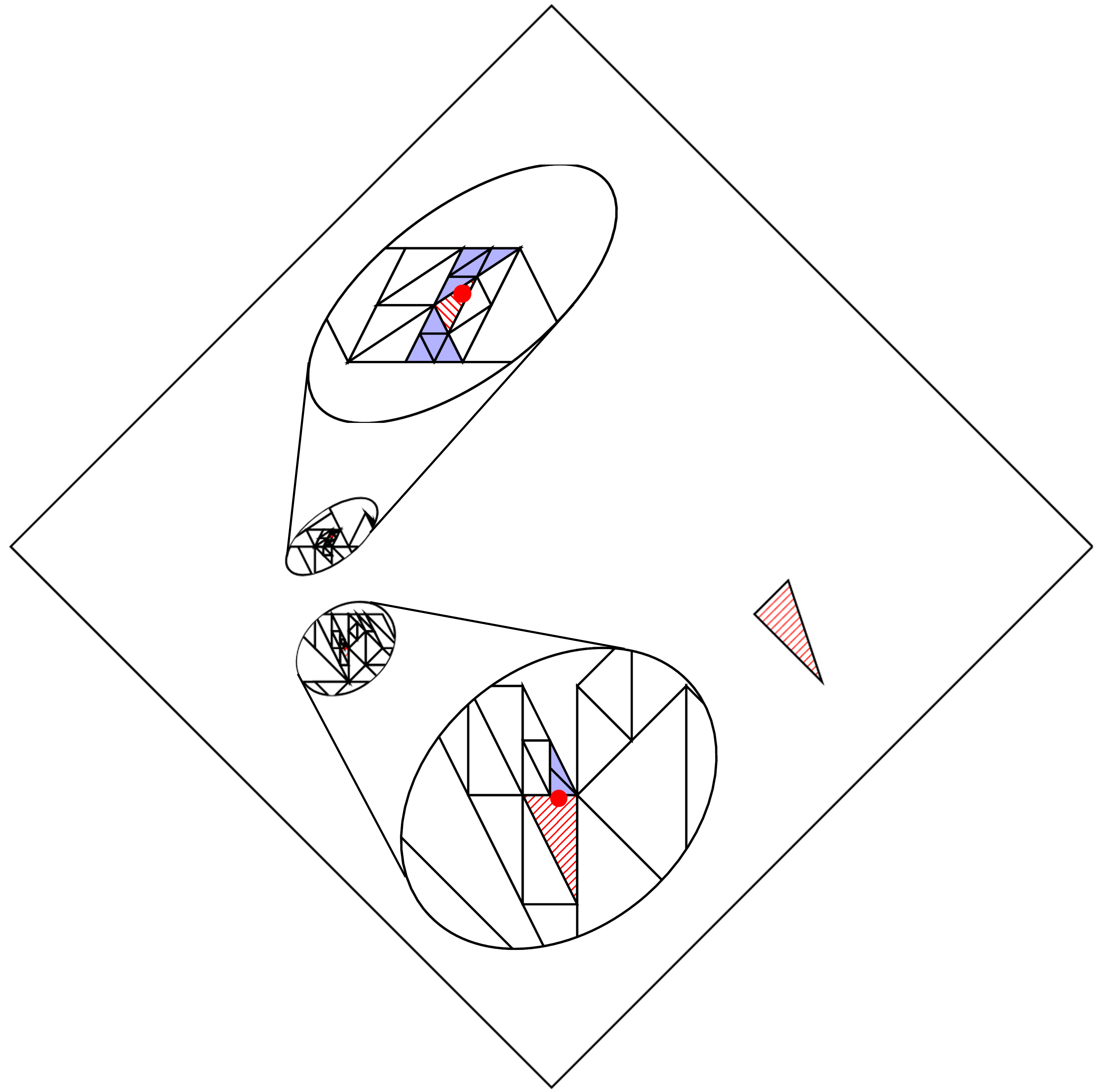}
    \caption{The triangulation with the zoomed-in regions for the two solutions}
    \label{fig:two-solutions2-result}
\end{figure}
\fi
We also chose the maximal and minimal area bounds to be $0.05$ and $6\cdot 10^{-5}$. The result is shown in \Cref{fig:two-solutions2-result} with the most interesting regions magnified.

From this result, we could see that there are three triangles labeled as ``accept'' found. Two of them contain the known solutions as shown in the zoomed-in region. For the other triangle, shaded in red, which is separated from the known solutions, we further approximate the solution with Newton's method by taking the center of mass of each triangle as the initial guess. At the end, we found one new solution at $\ww_3 = [\begin{matrix}
    .747677 & -.246394 & .616659
\end{matrix}]^T$, which is a third real solution to this scenario.

\subsection{Application to the IOD problem of a nearly circular orbit}

We also experimented with our method using data from a simulation of a nearly circular orbit with observers located on Earth.

\subsubsection*{Nearly circular Orbit}

We considered an example of a mission around Earth with a near-circular orbit, where the lines of sight are the lines connecting the observers on Earth and the corresponding orbit points of the simulated mission. This example orbit is inspired by the AQUA mission ~\cite{parkinson_aqua_2003}. 

For our experiment, we used the data from the work by Mancini et al. ~\cite{mancini_geometric_2023}. If we multiply the positions of observers by the Earth radii, we have the following position matrix with entries in kilometers.

\begin{equation*}
    \pp = \left[\begin{matrix}
     {1519}&{1143.89}&{1519}&-{3092}&-{818.1}\\
     -{4674}&-{6249.74}&-{4674}&{4873}&{4289}\\
     {4065.13}&{557.78}&-{4065.13}&-{2715.2}&{4649.09}
     \end{matrix}\right] 
\end{equation*}
and the corresponding direction matrix and the known solution are
\begin{equation*}
    \begin{aligned}
    \ww &= \left[\begin{matrix}
     -{.985693}\\
     -{.0898144}\\
     {.142629}
     \end{matrix}\right]\\
     \uu &= \left[\begin{matrix}
     -{.15563}&-{.379324}&-{.372229}&{.811237}&{.618925}\\
     {.460986}&-{.744255}&-{.661444}&{.574661}&{.775925}\\
     {.873654}&{.549725}&{.651105}&-{.107978}&-{.121953}
     \end{matrix}\right]
    \end{aligned}
\end{equation*}
\ifarxiv
\begin{figure}[H]
    \centering
    \includegraphics[width=0.5\linewidth]{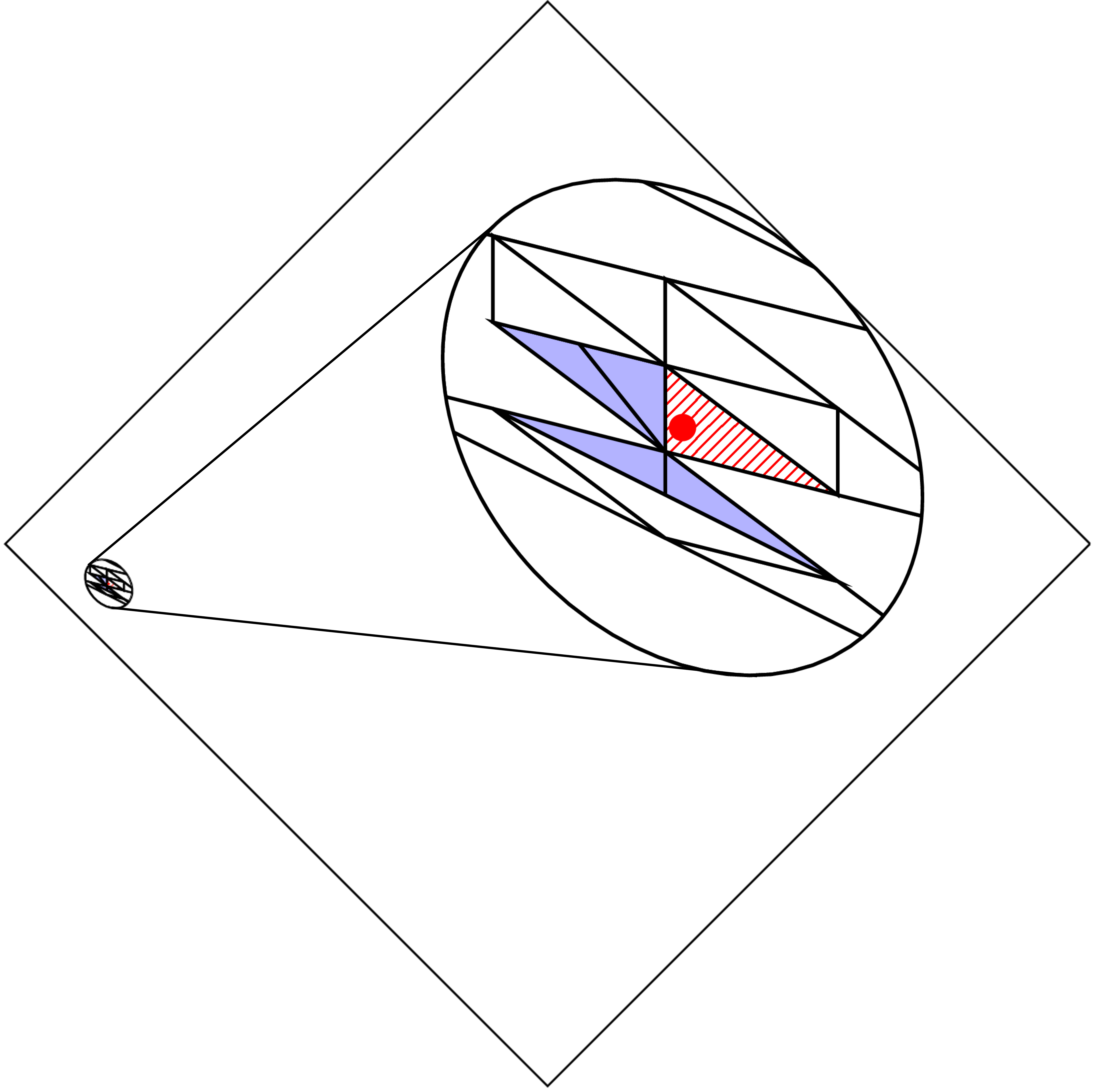}
    \caption{The triangulation from the nearly circular simulated data with the elliptical region around the true solution zoomed in}
    \label{figure:Nearly-Circular-result}
\end{figure}
\fi
Again, we used the oracle sequence $\oSequence$ in \Cref{eqn:OracleScheme}, and we used the same $\areaScalingThreshold$ and $\safetyCoefficient$ as the previous example but changed $\intersectionNormBound$ to $10000$.

We used the maximal and minimal area bounds of $0.05$ and $0.0003$ to start applying the oracles and terminate the subdivision. The triangulation and the zoomed-in elliptical region around the known solution are drawn in \Cref{figure:Nearly-Circular-result}. From the plot, we can see that we found one accepted triangle with some blue triangles around it, and we have the known solution inside the small accepted triangle we found. For this example, we also implemented both the certified oracles $\nonzeroCertifiedOracle$ and $\krawczykOralce$ to the rejected and accepted triangles. Practically, we first subdivide the accepted triangle into a set of very small triangles (with area less than $1\times 10^{19}$), then the Krawczyk certified oracle accepted one of these triangles certifiably, which means we proved the existence of a true solution in this accepted small triangle.

\subsection{Experiment statistics and computational efficiency}

\subsubsection{Experiment statistics}

From the experiments we described above, we compute the statistics as shown in \Cref{tab:stats}.

From the \Cref{tab:stats}, we could see that the ratios of the sum of the accepted and passed area with the rejected area are less than $0.05$ for all three cases. This means our projective plane subdivision method is able to refine the compact region of the four faces of the octahedron into a comparatively small triangular region containing the real solutions.

\begin{table}[ht]
\centering
\begin{tabular}{|c|c|c|c|}
\hline
\textbf{Data}                                                                 & \textbf{Single Observer} & \textbf{Two Solutions} & \textbf{Near Circular} \\ \hline
Area accepted                                                        & $0.003382$      & $0.014060$    & $0.000422$    \\ \hline
Area Passed                                                          & $0$             & $0.058355$    & $0.004757$   \\ \hline
\begin{tabular}[c]{@{}c@{}}Area rejected \\ by $\intersectionPointOracle$\end{tabular}  & $0.500670$      & $0.755763$    & $3.117353$   \\ \hline
\begin{tabular}[c]{@{}c@{}}Area rejected \\ by $\GDRejectingOracle$\end{tabular}     & $1.253791$      & $0.184262$    & $0.240926$   \\ \hline
\begin{tabular}[c]{@{}c@{}}Area rejected \\ by $\linearApproximationOracle$\end{tabular} & $1.706256$      & $2.451659$    & $0.100641$   \\ \hline
$\tfrac{\text{accept}+\text{pass}}{\text{reject}}$                                                & $0.000977$      & $0.021350$    & $0.001497$   \\ \hline
\end{tabular}
\caption{Statistics from the experiments}
\label{tab:stats}
\end{table}

\begin{table}[ht]
\centering
\begin{tabular}{|c|c|c|}
\hline
\textbf{Data}   & \textbf{\begin{tabular}[c]{@{}c@{}}Subdivision\end{tabular}} & \textbf{\begin{tabular}[c]{@{}c@{}}Homotopy Continuation\end{tabular}} \\ \hline
Single Observer & 6444                                                                            & 6015                                                                                                        \\ \hline
Two Solutions   & 28253                                                                          & 4118                                                                                                        \\ \hline
Nearly Circular    & 6583                                                                           & 4641                                                                                                        \\ \hline
\end{tabular}
\caption{Number of times the most expensive step in an algorithm is executed. (We believe that the cost of one step in homotopy continuation is \emph{higher by at least one order of magnitude} than one step of subdivision.)}
\label{tab:comparison}
\end{table}

\subsubsection{Comparison}

To compare our subdivision method to the homotopy continuation method for the dual formulation of the problem in~\cite{duff2022orbit} we report the number of calls to a subroutine which is a bottleneck for each method in
\Cref{tab:comparison}. 

The most expensive subroutine in the oracles used for subdivision is the evaluation of the Jacobian: this, in particular, relies on computing an inverse of a $5\times 5$ \emph{real} matrix.

In contrast, the predictor-corrector method of homotopy continuation in the referenced implementation takes from 5 to 7 linear algebra solver calls on the \emph{complex} instances of problems of size $7\times 7$. 

We should mention that our proof-of-concept implementation is carried out in the top-level interpreted language of Macaulay2~\cite{M2} and the core part of the referenced homotopy continuation routine is written in C++.
While for this reason, the direct timing comparison at this point is meaningless, the statistics we report suggest that a thorough implementation in a compiled language would be faster in practice on the examples of practical importance.

\section{Conclusion}
We produced a new direct (or primal) formulation for the geometric angles-only initial orbit determination problem and an algorithm based on a subdivision technique to locate its solutions. 

The experiments with our current implementation in Macaulay2 suggest that this approach is competitive with the existing method of Ref. \cite{mancini_geometric_2023}.

There is a range of flexibility in the choice of heuristic oracles that an engineer may use to take advantage of particular scenarios. An optimized implementation in a low-level compiled language is expected to be the fastest practical solver for this purely geometric problem.  

\ifarxiv
\section*{Acknowledgements}
\else 
\begin{credits}
\subsubsection{\ackname}
\fi
This work was partially supported by NSF DMS award 2001267 and AFOSR award FA95502310512. The authors would like to thank Dr. John Christian from the Daniel Guggenheim School of Aerospace Engineering of the Georgia Institute of Technology for his valuable guidance and insightful suggestions throughout the course of this research.
\ifarxiv
\else
\end{credits}
\fi
 This work was partially supported by NSF DMS award 2001267 and AFOSR award FA95502310512. The authors would like to thank Dr. John Christian from the Daniel Guggenheim School of Aerospace Engineering of the Georgia Institute of Technology for his valuable guidance and insightful suggestions throughout the course of this research.
\bibliographystyle{splncs04}
\bibliography{reference}
\end{document}